\documentclass[a4paper]{amsart}
\theoremstyle{plain}
  \newtheorem{thm}{Theorem}[section]
  \newtheorem{lem}[thm]{Lemma}
  \newtheorem{cor}[thm]{Corollary}
  \newtheorem{prop}[thm]{Proposition}
  \newtheorem{conj}[thm]{Conjecture}
  \newtheorem{claim}[thm]{Claim}

\theoremstyle{definition}

\theoremstyle{remark}
  \newtheorem{rem}[thm]{Remark}
  \newtheorem*{ack}{Acknowledgments}

\newcommand{\C}{\mathbb{C}}
\newcommand{\Vol}{\operatorname{Vol}}

\newcommand{\floor}[1]{\lfloor#1\rfloor}
\newcommand{\arccosh}{\operatorname{arccosh}}
\numberwithin{equation}{section}
\allowdisplaybreaks
\begin{document}
\title{Some limits of the colored Jones polynomials \\ of the figure-eight knot}
\author{Hitoshi Murakami}
\address{
Department of Mathematics,
Tokyo Institute of Technology,
Oh-okayama, Meguro, Tokyo 152-8551, Japan
}
\email{starshea@tky3.3web.ne.jp}
\date{\today}
\begin{abstract}
We will study the asymptotic behaviors of the colored Jones polynomials of
the figure-eight knot.
In particular we will show that for certain limits we obtain the volumes of
the cone manifolds with singularities along the knot.
\end{abstract}
\thanks{This research is partially supported by Grant-in-Aid for Scientific
Research (B)}
\maketitle
\section{Introduction}
Let $J_N(K;t)$ be the colored Jones polynomial of a knot $K$ associated with the
$N$-dimensional irreducible representation of the Lie algebra
$\mathfrak{sl}(2;\C)$,
normalized so that $J_N(\text{unknot};t)=1$.
If we evaluate it at the $N$-th root of unity
$\exp(2\pi\sqrt{-1}/N)$, then its absolute value coincides with
Kashaev's invariant
\cite{Kashaev:MODPLA95,
      Kashaev:LETMP97,
      Murakami/Murakami:ACTAM101}.
It was conjectured by R.~Kashaev \cite{Kashaev:LETMP97} that the growth rate of
his invariant for large $N$ determines the hyperbolic volume of a hyperbolic
knot.
His conjecture was generalized by J.~Murakami and the author to the following
Volume Conjecture:
\begin{conj}[Volume Conjecture]
Let $K$ be any knot.
Then
\begin{equation*}
  2\pi\lim_{N\to\infty}
  \frac{\log\left|J_N\left(K;\exp(2\pi\sqrt{-1}/N)\right)\right|}{N}
  =
  v_3\|S^3\setminus{K}\|,
\end{equation*}
where $v_3$ is the volume of the ideal regular hyperbolic tetrahedron
and $\|M\|$ is the Gromov norm (or simplicial volume).
In particular if $K$ is hyperbolic, then the right hand side coincides
with the volume of the knot complement $\Vol(S^3\setminus{K})$.
\end{conj}
See
\cite{
Kashaev/Tirkkonen:ZAPNS200,
Yokota:Murasugi70,
Yokota:volume00,
Yokota:Topology_Symposium2000,
Yokota:GTM02,
Yokota:INTIS03,
Murakami/Murakami/Okamoto/Takata/Yokota:EXPMA02,
Murakami:4_1,
Murakami:Nagoya99,
Murakami:SURIK00,
Murakami:SURIK02}
for related topics.
\par
In this paper we will study limits of the colored Jones polynomials
of the figure-eight knot evaluated at $\exp(2\pi r\sqrt{-1}/N)$ for a fixed
number $r$.
Put $\Lambda(z):=-\displaystyle\int_{0}^{z}\log|2\sin{x}|dx$,
the Lobachevski function,
and $\theta(r):=\arccos\left(\cos(2\pi r)-1/2\right)$ with
$0\le\arccos(x)\le\pi$.
\begin{thm}\label{thm}
Let $r$ be a real number satisfying ${5/6}<{r}<{7/6}$.
Then
\begin{equation*}
  2\pi\limsup_{N\to\infty}
  \frac{\log\left|J_N\left(E;\exp(2\pi r\sqrt{-1}/N)\right)\right|}{N}
  =
  \frac{2\Lambda\bigl(\pi r+\theta(r)/2\bigr)
       -2\Lambda\bigl(\pi r-\theta(r)/2\bigr)}{r}.
\end{equation*}
\par
Moreover if $r$ is irrational or $r=1$, then
\begin{equation*}
  2\pi\lim_{N\to\infty}
  \frac{\log\left|J_N\left(E;\exp(2\pi r\sqrt{-1}/N)\right)\right|}{N}
  =
  \frac{2\Lambda\bigl(\pi r+\theta(r)/2\bigr)
       -2\Lambda\bigl(\pi r-\theta(r)/2\bigr)}{r},
\end{equation*}
and if $r\ne1$ and rational, then
\begin{equation*}
  2\pi\liminf_{N\to\infty}
  \frac{\log\left|J_N\left(E;\exp(2\pi r\sqrt{-1}/N)\right)\right|}{N}
  =0
\end{equation*}
\end{thm}
\begin{rem}
The case where $r=1$ is due to R.~Kashaev \cite{Kashaev:LETMP97}
and T.~Ekholm \cite{Murakami:4_1}, and the following proof is similar to
Ekholm's.
\end{rem}
\begin{rem}
The value
$2\Lambda\bigl(\pi r+\theta(r)/2\bigr)-2\Lambda\bigl(\pi r-\theta(r)/2\bigr)$
coincides with the volume of the cone-manifold with underlying space $S^3$,
the singularity the figure-eight knot, and the cone-angle $2\pi|1-r|$
\cite{Vesnin/Mednykh:1995,Gukov:03},
which was informed by S. Gukov, A. Mednykh, and A. Vesnin.
\end{rem}
\begin{rem}
Some results in \cite{Murakami:SURIK02} were erroneously stated.
The author did not consider the case where $r$ is rational.
\end{rem}
We also calculate the limits for some other cases.
\section{Preliminaries}
Let $E$ denote the figure-eight knot $4_1$.
Due to K.~Habiro and T.~Le, the following formula is known.
\begin{equation}\label{eq:fig8}
  J_N(E;t)
  =
  \sum_{k=0}^{N-1}\prod_{j=1}^{k}
  \left(
    t^{(N+j)/2}-t^{-(N+j)/2}
  \right)
  \left(
    t^{(N-j)/2}-t^{-(N-j)/2}
  \right).
\end{equation}
Put $t=\exp(2\pi r\sqrt{-1}/N)$.
Since
\begin{align*}
  t^{(N+j)/2}-t^{-(N+j)/2}&=2\sqrt{-1}\sin(\pi r(N+j)/N)
  \\
  \intertext{and}
  t^{(N-j)/2}-t^{-(N-j)/2}&=2\sqrt{-1}\sin(\pi r(N-j)/N),
\end{align*}
we have
\begin{multline*}
  \left(
    t^{(N+j)/2}-t^{-(N+j)/2}
  \right)
  \left(
    t^{(N-j)/2}-t^{-(N-j)/2}
  \right)
  \\
  =
  4\sin(\pi rj/N+\pi r)\sin(\pi rj/N-\pi r).
\end{multline*}
Put
\begin{align*}
g(j):&=4\sin(\pi rj/N+\pi r)\sin(\pi rj/N-\pi r)
\\
     &=2\cos(2\pi r)-2\cos(2\pi rj/N)
\\
\intertext{and}
f(k):&=\prod_{j=1}^{k}g(j)
\end{align*}
so that
$J_N(E;\exp(2\pi r\sqrt{-1}/N))=\sum_{k=0}^{N-1}f(k)$.
\section{The case where $r$ is irrational and ${5/6}<r<{1}$ or $r=1$.}
In this section we prove Theorem~\ref{thm} in the case where $5/6<r<1$ and $r$ is irrational, or $r=1$.
\par
We put
\begin{equation*}
  B:=\frac{N(1-r)}{r},\quad
  C:=\frac{N\theta(r)}{2\pi r},\quad
  \text{and}\quad
  D:=\frac{N(2\pi-\theta(r))}{2\pi r}
\end{equation*}
with $\theta(r):=\arccos\left(\cos(2\pi r)-1/2\right)$.
Here $\arccos$ takes its value between $0$ and $\pi$.
Note that $0\le B<C<D<1$ ($B=0$ only if $r=1$),
that $g(B)=0$ and $g(C)=g(D)=1$,
and that $B$ is not an integer for any $N$ since $r$ is irrational.
Therefore $g(j)\ne0$ for any integer $0<j<N$ in this case.
\par
Since we have
\begin{equation*}
\begin{array}{rccr}
  -1<g(j)<0 & \text{for} & 0<j<B, &
  \\
  0<g(j)<1  & \text{for} & B<j<C, &
  \\
  g(j)>1    & \text{for} & C<j<D, & \text{and}
  \\
  0<g(j)<1  & \text{for} & D<j<N, &
\end{array}
\end{equation*}
\par\bigskip
\begin{center}
\begin{tabular}{|c||c|c|c|c|c|c|c|c|c|}
  \hline
  \multicolumn{10}{|c|}{Table of $g(j)$ regarding $j$ as a continuous parameter}
  \\
  \hline\hline
  $j$&$0$&$\cdots$&$B$&$\cdots$&$C$&$\cdots$&$D$&$\cdots$&$N$
  \\
  \hline
  $g(j)$&$-1<g(0)<0$&$\nearrow$&$0$&$\nearrow$&$1$&$>1$&$1$&$\searrow$&$0$
  \\
  \hline
\end{tabular}
\end{center}
\par\bigskip
we see
\begin{enumerate}
\item
  If $j<B$ then the signs of $f(j)$ alternate, that is,
  $f(j-1)f(j)<0$, and if $j>B$ then the signs of $f(j)$ are constant, and
\item
  $1=|f(0)|>|f(1)|>\dots>|f(\floor{C}-1)|$
  and
  $|f(\floor{C})|<\dots<|f(\floor{D})|$, where $\floor{x}$ is the greatest
  integer that does not exceed $x$.
\end{enumerate}
Put $F_N:=|f(\floor{D})|$, which is the maximum of
$\{|f(j)|\}$ for $0\le j \le N-1$.
\par
We will show the following inequality.
\begin{claim}\label{claim:inequalities5/6<r<1}
\begin{equation*}
  F_N-1 \le \left|J_N(E;\exp(2\pi r\sqrt{-1}/N))\right| \le NF_N
\end{equation*}
\end{claim}
\begin{proof}
The last inequality follows since $|f(j)|\le F_N$ and
$J_N(E;\exp(2\pi r\sqrt{-1}/N))=\sum_{j=0}^{N-1}f(j)$.
\par
To prove the first inequality, first we consider the case where $\floor{B}$ is
odd.
In this case since $f(0)=1$, $f(2j-1)+f(2j)<0$ for $2j<\floor{B}$,
and $f(j)<0$ for
$j\ge\floor{B}$, we have
\begin{align*}
  &\left|J_N(E;\exp(2\pi r\sqrt{-1}/N))\right|
  \\
  &\quad=
  \bigl|
  f(0)+
  \{f(1)+f(2)\}+\{f(3)+f(4)\}+\dots+\{f(\floor{B}-2)+f(\floor{B}-1)\}
  \\
  &\quad\quad
  +f(\floor{B})+f(\floor{B}+1)+\dots+f(N-1)
  \bigr|
  \\
  &\quad=
  -1+
  |f(1)+f(2)|+|f(3)+f(4)|+\dots+|f(\floor{B}-2)+f(\floor{B}-1)|
  \\
  &\quad\quad
  +|f(\floor{B})|+|f(\floor{B}+1)|+\dots+|f(\floor{D})|+\dots+|f(N-1)|
  \\
  &\quad>
  F_N-1
\end{align*}
and the first equality follows.
\par
Next we consider the case where $\floor{B}$ is even.
In this case since $f(2j)+f(2j+1)>0$ for $2j+1<\floor{B}$,
and $f(j)>0$ for $j\ge\floor{B}$, we have
\begin{align*}
  &\left|J_N(E;\exp(2\pi r\sqrt{-1}/N))\right|
  \\
  &\quad=
  \bigl|
  \{f(0)+f(1)\}+\{f(2)+f(3)\}+\dots+\{f(\floor{B}-2)+f(\floor{B}-1)\}
  \\
  &\quad\quad
  +f(\floor{B})+f(\floor{B}+1)+\dots+f(N-1)
  \bigr|
  \\
  &\quad=
  |f(0)+f(1)|+|f(2)+f(3)|+\dots+|f(\floor{B}-2)+f(\floor{B}-1)|
  \\
  &\quad\quad
  +|f(\floor{B})|+|f(\floor{B}+1)|+\dots+|f(\floor{D})|+\dots+|f(N-1)|
  \\
  &\quad>
  F_N.
\end{align*}
and the first equality also follows.
\end{proof}
Now we study the asymptotic behavior of $F_N$.
\begin{claim}
\begin{equation*}
  \lim_{N\to\infty}\frac{\log{F_N}}{N}
  =
  \frac{1}{\pi r}
  \left\{
    \Lambda\left(\pi r+\theta(r)/2\right)-\Lambda\left(\pi r-\theta(r)/2\right)
  \right\}.
\end{equation*}
\end{claim}
\begin{proof}
Since
\begin{align*}
  F_N&=|f(\floor{D})|
  =
  \prod_{j=1}^{\floor{D}}|g(j)|
  \\
  &=
  \prod_{j=1}^{\floor{D}}
  \left|2\sin(\pi rj/N+\pi r)\right|
  \left|2\sin(\pi rj/N-\pi r)\right|,
\end{align*}
we have
\begin{align*}
  &\lim_{N\to\infty}\frac{\log{F_N}}{N}
  \\
  &\quad=
  \lim_{N\to\infty}\frac{1}{N}
  \sum_{j=1}^{\floor{D}}
  \left\{
     \log\left|2\sin(\pi rj/N+\pi r)\right|
    +\log\left|2\sin(\pi rj/N-\pi r)\right|
  \right\}
  \\
  &\quad=
   \frac{1}{\pi r}
   \int_{ \pi r}^{\pi-\theta(r)/2+\pi r}\log|2\sin{x}|d\,x
  +\frac{1}{\pi r}
   \int_{-\pi r}^{\pi-\theta(r)/2-\pi r}\log|2\sin{x}|d\,x
  \\
  &\quad=
  \frac{1}{\pi r}
  \left\{
    -\Lambda(\pi-\theta(r)/2+\pi r)+\Lambda( \pi r)
    -\Lambda(\pi-\theta(r)/2-\pi r)+\Lambda(-\pi r)
  \right\}
  \\
  &\quad=
  \frac{1}{\pi r}
  \left\{
    \Lambda(\pi r+\theta(r)/2)-\Lambda(\pi r-\theta(r)/2)
  \right\}.
\end{align*}
Here we use the formulas $-\Lambda(-x)=\Lambda(x+\pi)=\Lambda(x)$.
\end{proof}
Since $\Lambda(\pi r+\theta(r)/2)-\Lambda(\pi r-\theta(r)/2)>0$ (see Appendix),
$F_N$ grows exponentially and so
$\lim_{N\to\infty}\log(F_N-1)/N=\lim_{N\to\infty}\log{F_N}/N$.
\par
Now since $\lim_{N\to\infty}\log{N}/N=0$, we have from
Claim~\ref{claim:inequalities5/6<r<1}
\begin{align*}
  \lim_{N\to\infty}\frac{\log|J_N(E;\exp(2\pi r\sqrt{-1}/N))|}{N}
  &=
  \lim_{N\to\infty}\frac{\log{F_N}}{N}
  \\
  &=
  \frac{1}{\pi r}
  \left\{
    \Lambda(\pi r+\theta(r)/2)-\Lambda(\pi r-\theta(r)/2)
  \right\}.
\end{align*}
This completes proof of the theorem in the case where $r$ is irrational with
${5/6}<r<{1}$ or $r=1$.
\section{The case where $r$ is irrational and $1<r<7/6$}
We put
\begin{equation*}
  B:=\frac{N(r-1)}{r},\quad
  C:=\frac{N\theta(r)}{2\pi r},\quad
  D:=\frac{N(2\pi-\theta(r))}{2\pi r},\quad
  \text{and}\quad
  B':=\frac{N(2-r)}{r}
\end{equation*}
Note that $0<B<C<D<1$,
that $g(B)=g(B')=0$ and $g(C)=g(D)=1$,
and that $B$ and $B'$ are not integers for any $N$ since $r$ is irrational.
Therefore $g(j)\ne0$ for any integer $0<j<N$ in this case.
\par
Since we have
\begin{equation*}
\begin{array}{rccr}
  -1<g(j)<0 & \text{for} & 0<j<B,  &
  \\
  0<g(j)<1  & \text{for} & B<j<C,  &
  \\
  g(j)>1    & \text{for} & C<j<D,  &
  \\
  0<g(j)<1  & \text{for} & D<j<B', & \text{and}
  \\
  -1<g(j)<0 & \text{for} & B'<j<N, &
\end{array}
\end{equation*}
the proof is similar to the case where $5/6<r<1$.
We put $F_N:=\left|f(\floor{D})\right|$.
\begin{center}
\begin{tabular}{|c||c|c|c|c|c|c|c|c|c|c|c|}
  \hline
  \multicolumn{12}{|c|}{Table of $g(j)$ regarding $j$ as a continuous parameter}
  \\
  \hline\hline
  $j$&$0$&$\cdots$&$B$&$\cdots$&$C$&$\cdots$&$D$&$\cdots$&$B'$&$\cdots$&$N$
  \\
  \hline
  $g(j)$&$-1<g(0)<0$&$\nearrow$&$0$&$\nearrow$&$1$&$g(j)>1$&$1$&$\searrow$&$0$&$-1<g(j)<0$&$0$
\\
\hline
\end{tabular}
\end{center}
If $\floor{B}$ is odd, we have
\begin{align*}
  &\left|J_N(E;\exp(2\pi r\sqrt{-1}/N))\right|
  \\
  &\quad=
  \bigl|
  f(0)+
  \{f(1)+f(2)\}+\{f(3)+f(4)\}+\dots+\{f(\floor{B}-2)+f(\floor{B}-1)\}
  \\
  &\quad\quad
  +f(\floor{B})+f(\floor{B}+1)+\dots+f(\floor{B'}-1)
  \\
  &\quad\quad
  +f(\floor{B'})+f(\floor{B'}+1)+\dots+f(N-1)
  \bigr|
  \\
  &\quad=
  -1+
  |f(1)+f(2)|+|f(3)+f(4)|+\dots+|f(\floor{B}-2)+f(\floor{B}-1)|
  \\
  &\quad\quad
  +|f(\floor{B})|+|f(\floor{B}+1)|+\dots+|f(\floor{D})|+\dots+|f(\floor{B'}-1)|
  \\
  &\quad\quad
  +
  \begin{cases}
    |f(\floor{B'})+f(\floor{B'}+1)|+\dots+|f(N-2)-f(N-1)|
    \\
    \hfill\text{if $N-\floor{B'}$ is even}
    \\
    |f(\floor{B'})+f(\floor{B'}+1)|+\dots+|f(N-3)-f(N-2)|+|f(N-1)|
    \\
    \hfill\text{if $N-\floor{B'}$ is odd}
  \end{cases}
  \\
  &\quad>
  F_N-1
\end{align*}
\par
If $\floor{B}$ is even,
\begin{align*}
  &\left|J_N(E;\exp(2\pi r\sqrt{-1}/N))\right|
  \\
  &\quad=
  \bigl|
  \{f(0)+f(1)\}+\{f(2)+f(3)\}+\dots+\{f(\floor{B}-2)+f(\floor{B}-1)\}
  \\
  &\quad\quad
  +f(\floor{B})+f(\floor{B}+1)+\dots+f(\floor{B'}-1)
  \\
  &\quad\quad
  +f(\floor{B'})+|f(\floor{B'}+1)|+\dots+|f(N-1)|
  \bigr|
  \\
  &\quad=
  |f(0)+f(1)|+|f(2)+f(3)|+\dots+|f(\floor{B}-2)+f(\floor{B}-1)|
  \\
  &\quad\quad
  +|f(\floor{B})|+|f(\floor{B}+1)|+\dots+|f(\floor{D})|+\dots+|f(\floor{B'}-1)|
  \\
  &\quad\quad
  +
  \begin{cases}
    |f(\floor{B'})+f(\floor{B'}+1)|+\dots+|f(N-2)-f(N-1)|
    \\
    \hfill\text{if $N-\floor{B'}$ is even}
    \\
    |f(\floor{B'})+f(\floor{B'}+1)|+\dots+|f(N-3)-f(N-2)|+|f(N-1)|
    \\
    \hfill\text{if $N-\floor{B'}$ is odd}
  \end{cases}
  \\
  &\quad>
  F_N.
\end{align*}
\par
So we have $F_N-1<|J_N(E;\exp(2\pi r\sqrt{-1}/N))|<NF_N$, and the result follows, completing
the proof of Theorem~\ref{thm}.
\par
As a corollary we have
\begin{cor}\label{cor:limit}
\begin{equation*}
  \lim_{r\to1}
  \left\{
    \limsup_{N\to\infty}
    \frac{\log\left|J_N(E;\exp(2\pi r\sqrt{-1}))\right|}{N}
  \right\}
  =
    \limsup_{N\to\infty}
    \frac{\log\left|J_N(E;\exp(2\pi\sqrt{-1}))\right|}{N}
\end{equation*}
\end{cor}
\section{The case where $r$ is rational and $|1-r|<{1/6}$.}
Put $r:=q/p$ with coprime integers $p$ and $q$.
\par
Suppose first that $N$ is not a multiple of $q$.
Then $B$ and $B'$ are not integers and $g(j)\ne0$ for any $j$.
So in this case a similar argument for the case where $r$ is irrational applies
and we have
\begin{equation*}
  \frac{\log|J_{N'}(E;\exp(2\pi r\sqrt{-1}/N))|}{N'}
  \to
  \frac{1}{\pi r}
  \left\{
    \Lambda(\pi r+\theta(r)/2)-\Lambda(\pi r-\theta(r)/2)
  \right\}
\end{equation*}
for the subsequence $\{N'\}$ of natural numbers which are not multiples of $q$.
\par
Next suppose that $N$ is a multiple of $q$, say $N=nq$.
Then $B$ is an integer since $B=n(p-q)$ ($n(q-p)$ respectively) if $5/6<r<1$ ($1<r<7/6$ respectively).
Therefore $-1<g(j)<0$ for $0<j<B$ and $g(j)=0$ for $j\ge B$ and so
$|f(j)|<1$ for $0<j<B$ and $f(j)=0$ for $j\ge B$.
Thus we have
\begin{multline*}
  \left|J_N(E;\exp(2\pi r\sqrt{-1}/N))\right|
  =
  \left|\sum_{k=0}^{N-1}f(k)\right|
  =
  \left|\sum_{k=0}^{B-1}f(k)\right|
  \le
  \sum_{k=0}^{B-1}|f(k)|
  \\
  <B=N|1-1/r|.
\end{multline*}
Moreover since
\begin{align*}
  \left|\sum_{k=0}^{B-1}f(k)\right|
  &=
  \begin{cases}
    \{1+f(1)\}+\{f(2)+f(3)\}+\dots+\{f(\floor{B}-2)+f(\floor{B}-1)\}
    \\
    \hfill\text{if $\floor{B}$ is even}
    \\
    \{1+f(1)\}+\dots+\{f(\floor{B}-3)+f(\floor{B}-2)\}+f(\floor{B}-1)
    \\
    \hfill\text{if $\floor{B}$ is odd}
  \end{cases}
  \\
  &>
  1+f(1)
  =
  1+2\cos(2\pi r)-2\cos(2\pi r/N)
  \\
  &>2-2\cos(2\pi r/N).
\end{align*}
Therefore
\begin{multline*}
  |J_N(E;\exp(2\pi r\sqrt{-1}/N))|
  >
  2-2\cos(2\pi r/N)
  >
  2-2\left(1-\frac{(2\pi r/N)^6}{6!}\right)
  \\
  =
  (2\pi r)^6/360N^6.
\end{multline*}
So we finally have
\begin{equation*}
  \frac{\log\left((2\pi r)^6/360N^6\right)}{N}
  <
  \frac{\log|J_N(E;\exp(2\pi r\sqrt{-1}/N))|}{N}
  <
  \frac{\log(N|1-1/r|)}{N}
\end{equation*}
and so $\log|J_N(E;\exp(2\pi r\sqrt{-1}/N))|/N$ can be arbitrarily small in this case.
Therefore
\begin{equation*}
  \frac{\log|J_{N''}(E;\exp(2\pi r\sqrt{-1}/N))|}{N''}
  \to
  0
\end{equation*}
for the subsequence $\{N''\}$ of multiples of $q$.
\par
The proof in this case is complete.
\section{The case where $0\le{r}<{1/6}$.}
In this section we will show the following proposition.
\begin{prop}\label{prop:r<1/6}
For any $r$ with $0\le r<1/6$, we have
\begin{equation*}
  \lim_{N\to\infty}\frac{\log|J_N(E;\exp(2\pi r\sqrt{-1}/N))|}{N}=0.
\end{equation*}
\end{prop}
\par
\begin{proof}
Since $J_N(E;1)=1$, we have the equality when $r=0$.
\par
If $0<r<{1/6}$, then $-1<g(j)<0$ for $j=1,2,\dots,N-1$ and so
$|f(k)|<1$ for $k=1,2,\dots,N-1$.
Therefore
\begin{equation*}
  \left|J_N\left(E;\exp(2\pi r\sqrt{-1}/N)\right)\right|
  =
  \left|\sum_{k=0}^{N-1}f(k)\right|
  <N.
\end{equation*}
On the other hand, since $f(k)$ is positive (negative, respectively)
if $k$ is even (odd, respectively) and $|f(k)|$ is decreasing, we have
\begin{align*}
  &\left|\sum_{k=0}^{N-1}f(k)\right|
  \\
  &\quad=
  \begin{cases}
    \{1+f(1)\}+\{f(2)+f(3)\}+\dots+\{f(N-2)+f(N-1)\}
    \\
    \hfill\text{if $N$ is even}
    \\
    \{1+f(1)\}+\{f(2)+f(3)\}+\dots+\{f(N-3)+f(N-2)\}+f(N-1)
    \\
    \hfill\text{if $N$ is odd}
  \end{cases}
  \\
  &\quad>
  1+f(1)
  =1+2\cos(2\pi r)-2\cos(2\pi r/N)
  >
  (2\pi r)^6/360N^6.
\end{align*}
Therefore we have
\begin{equation*}
  \frac{\log\left((2\pi r)^6/360N^6\right)}{N}
  <
  \frac{\log\left(J_N\left(E;\exp(2\pi r\sqrt{-1}/N)\right)\right)}{N}
  <
  \frac{\log N}{N}
\end{equation*}
and the proposition follows.
\end{proof}
\section{The case where $r$ is irrational and $1/6\le|1-r|<1/4$.}
In this section we will prove a result similar to Theorem\ref{thm} for
irrational $r$ with $1/6\le|1-r|<1/4$.
\begin{prop}\label{prop:1/6<|1-r|<1/4}
If $r$ is irrational and $1/6\le|1-r|<1/4$, then we have
\begin{equation*}
  \limsup_{N\to\infty}
  2\pi\frac{\log\left|J_N\left(E;\exp(2\pi r\sqrt{-1}/N)\right)\right|}{N}
  =
  \frac{2\Lambda\bigl(\pi r+\theta(r)/2\bigr)
       -2\Lambda\bigl(\pi r-\theta(r)/2\bigr)}{r}.
\end{equation*}
\end{prop}
First we prove the proposition in the case where $3/4<r\le5/6$.
\begin{proof}[Proof when $3/4<r\le5/6$.]
We put
\begin{equation*}
  A:=\frac{N\varphi(r)}{2\pi r},\quad
  B:=\frac{N(1-r)}{r},\quad
  C:=\frac{N\theta(r)}{2\pi r},\quad
  \text{and}\quad
  D:=\frac{N(2\pi-\theta(r))}{2\pi r},
\end{equation*}
with $\varphi(r):=\arccos\left(\cos(2\pi r)+1/2\right)$ and
$\theta(r):=\arccos\left(\cos(2\pi r)-1/2\right)$.
Note that $g(A)=-1$, $g(B)=0$, and $g(C)=g(D)=1$.
\par
Since
\begin{equation*}
\begin{array}{rccr}
  g(j)<-1   & \text{for} & j<A,   &
  \\
  -1<g(j)<0 & \text{for} & A<j<B, &
  \\
  0<g(j)<1  & \text{for} & B<j<C, &
  \\
  g(j)>1    & \text{for} & C<j<D, & \text{and}
  \\
  0<g(j)<1 & \text{for} & D<j,   &
\end{array}
\end{equation*}
\begin{center}
\begin{tabular}{|c||c|c|c|c|c|c|c|c|c|c|c|c|}
  \hline
  \multicolumn{12}{|c|}{Table of $g(j)$ regarding $j$ as a continuous parameter}
  \\
  \hline\hline
  $j$&$0$&$\cdots$&$A$&$\cdots$&$B$&$\cdots$&$C$&$\cdots$&$D$&$\cdots$&$N$
  \\
  \hline
  $g(j)$&$g(0)<-1$&$\nearrow$&$-1$&$\nearrow$&$0$&$\nearrow$&$1$&
  $g(j)>1$&$1$&$\searrow$&$0$
\\
\hline
\end{tabular}
\end{center}
we have
\begin{itemize}
\item
$1=|f(0)|<|f(1)|<\dots<|f(\floor{A}-1)|<|f(\floor{A})|$,
\item
$|f(\floor{A})|>|f(\floor{A}+1)|>\dots>|f(\floor{B}-1)|>|f(\floor{B})|
>|f(\floor{B}+1)|>\dots>|f(\floor{C}-1)|>|f(\floor{C})|$,
\item
$|f(\floor{C})|<|f(\floor{C}+1)|<\dots<|f(\floor{D}-1)|<|f(\floor{D})|$, and
\item
$|f(\floor{D})|>|f(\floor{D}+1)|>\dots>|f(N-2)|>|f(N-1)|$.
\end{itemize}
Therefore $|f(k)|$ takes its `local maxima' at $k=\floor{A}$ and $\floor{D}$.
\par
Now we have
\begin{align*}
  &
  \lim_{N\to\infty}\frac{\log\left|f(\floor{A})\right|}{N}
  \\
  =&
  \lim_{N\to\infty}
  \frac{1}{N}
  \sum_{j=1}^{\floor{A}}
  \left\{\log|2\sin(\pi rj/N+\pi r)|+\log|2\sin(\pi rj/N-\pi r)|\right\}
  \\
  =&
  \frac{1}{\pi r}\int_{ \pi r}^{\varphi(r)/2+\pi r}\log|2\sin y|\,dy
  +
  \frac{1}{\pi r}\int_{-\pi r}^{\varphi(r)/2-\pi r}\log|2\sin y|\,dy
  \\
  =&
  \frac{\Lambda(\pi r-\varphi(r)/2)-\Lambda(\pi r+\varphi(r)/2)}{\pi r}
\end{align*}
and similarly
\begin{align*}
  &
  \lim_{N\to\infty}\frac{\log\left|f(\floor{D})\right|}{N}
  \\
  =&
  \frac{1}{\pi r}\int_{ \pi r}^{\pi-\theta(r)/2+\pi r}\log|2\sin y|\,dy
  +
  \frac{1}{\pi r}\int_{-\pi r}^{\pi-\theta(r)/2-\pi r}\log|2\sin y|\,dy
  \\
  =&
  \frac{\Lambda(\pi r+\theta(r)/2)-\Lambda(\pi r-\theta(r)/2)}{\pi r}.
\end{align*}
\par
Since
\begin{multline}\label{eq:delta}
  \delta:=
  \pi r
  \left(
    \lim_{N\to\infty}\frac{\log|f(\floor{D})|}{N}
    -
    \lim_{N\to\infty}\frac{\log|f(\floor{A})|}{N}
  \right)
  \\
  =
  \Lambda(\pi r+\theta(r)/2)-\Lambda(\pi r-\theta(r)/2)
  +
  \Lambda(\pi r+\varphi(r)/2)-\Lambda(\pi r-\varphi(r)/2)
  >0,
\end{multline}
(see Appendix) we see that $|f(\floor{D})|$ is the maximum, that is,
$|f(k)|\le|f(\floor{D})|$ for sufficient large $N$.
Therefore $|J_N(E;\exp(2\pi r\sqrt{-1}/N))|<N|f(\floor{D})|$.
\par
On the other hand, since the sign of $f(k)$ does not change for $k\ge\floor{B}$
and $|f(k)|<|f(\floor{A})|$ for $k\le\floor{C}$, we have
\begin{multline*}
  |J_N(E;\exp(2\pi r\sqrt{-1}/N))|
  \\
  =
  \left|\sum_{k=0}^{N-1}f(k)\right|
  >
  \left|f(\floor{D})-\sum_{k=0}^{\floor{B-1}}f(k)\right|
  >
  \left|f(\floor{D})-Nf(\floor{A})\right|.
\end{multline*}
Since
\begin{equation*}
  \frac{|f(\floor{A})|}{|f(\floor{D})|}
  \underset{N\to\infty}{\sim}
  \exp(-\delta N)
\end{equation*}
from \eqref{eq:delta}, we have
\begin{align*}
  |J_N(E;\exp(2\pi r\sqrt{-1}/N))|
  &>
  |f(\floor{D})|\left(1-N\frac{f(\floor{A})}{f(\floor{D})}\right)
  \\
  &\underset{N\to\infty}{\sim}
  |f(\floor{D})|\{1-N\exp(-\delta N)\}
  \\
  &\xrightarrow{N\to\infty}
  |f(\floor{D})|.
\end{align*}
So we finally have
\begin{multline*}
  \lim_{N\to\infty}\frac{\log|J_N(E;\exp(2\pi r\sqrt{-1}/N))|}{N}
  =
  \lim_{N\to\infty}\frac{\log|f(\floor{D})|}{N}
  \\
  =
  \frac{\Lambda(\pi r+\theta(r)/2)-\Lambda(\pi r-\theta(r)/2)}{\pi r}
\end{multline*}
as required.
\end{proof}
Next we will consider the case where $7/6\le r<5/4$.
\begin{proof}[Proof when $7/6\le r<5/4$.]
We put
\begin{multline*}
  A:=\frac{N\varphi(r)}{2\pi r},\quad
  B:=\frac{N(r-1)}{r},\quad
  C:=\frac{N\theta(r)}{2\pi r},\quad
  D:=\frac{N(2\pi-\theta(r))}{2\pi r},\quad
  \\
  B':=\frac{N(2-r)}{r},\quad
  A':=\frac{N(2\pi-\varphi(r))}{r},\quad
  \text{and}\quad
  A'':=\frac{N(2\pi+\varphi(r))}{r}
\end{multline*}
with $\varphi(r):=\arccos\left(\cos(2\pi r)+1/2\right)$ and
Note that $g(A)=g(A')=g(A'')=-1$, $g(B)=g(B')=0$, and $g(C)=g(D)=1$.
\par
Since
\begin{equation*}
\begin{array}{rccr}
  g(j)<-1   & \text{for} & j<A,   &
  \\
  -1<g(j)<0 & \text{for} & A<j<B, &
  \\
  0<g(j)<1  & \text{for} & B<j<C, &
  \\
  g(j)>1    & \text{for} & C<j<D, &
  \\
  0<g(j)<1  & \text{for} & D<j<B', &
  \\
  -1<g(j)<0 & \text{for} & B'<j<A', &
  \\
  g(j)<-1   & \text{for} & A'<j<A'', & \text{and}
  \\
  -1<g(j)<0 & \text{for} & j>A'',
\end{array}
\end{equation*}
$|f(k)|$ has three `local maxima' at $k=\floor{A}$, $\floor{D}$, and
$\floor{A''}$.
\par
But since
\begin{align*}
  &
  \lim_{N\to\infty}\frac{\log\left|f(\floor{A''})\right|}{N}
  \\
  =&
  \lim_{N\to\infty}
  \frac{1}{N}
  \sum_{j=1}^{\floor{A''}}
  \left\{\log|2\sin(\pi rj/N+\pi r)|+\log|2\sin(\pi rj/N-\pi r)|\right\}
  \\
  =&
  \frac{1}{\pi r}\int_{ \pi r}^{\pi+\varphi(r)/2+\pi r}\log|2\sin y|\,dy
  +
  \frac{1}{\pi r}\int_{-\pi r}^{\pi+\varphi(r)/2-\pi r}\log|2\sin y|\,dy
  \\
  =&
  \frac{\Lambda(\pi r-\varphi(r)/2)-\Lambda(\pi r+\varphi(r)/2)}{\pi r}
  <0,
\end{align*}
$|f(\floor{D})|$ is the unique maximum.
Therefore the argument for the case $3/4<r\le5/6$ can be applied and the proof is complete.
\end{proof}
\section{The case where $r$ is purely imaginary.}
In this section we consider the case where $t=\exp(2\pi r/N)$.
\begin{thm}
If $2\pi|r|>\arccosh(3/2)$, then
\begin{equation*}
  2\pi|r|\lim_{N\to\infty}
  \frac{\log J_N(E;\exp(2\pi r/N))}{N}
  =
  2\Gamma(\pi|r|+\varphi(|r|)/2)-2\Gamma(\pi|r|-\varphi(|r|)/2),
\end{equation*}
where
$\Gamma(z):=\int_{0}^{z}\log 2\sinh x d\,x$
and
$\varphi(r):=\arccosh(\cosh(2\pi r)-1/2)$.
If $2\pi|r|\le\arccosh(3/2)$, then
\begin{equation*}
  \lim_{N\to\infty}
  \frac{\log J_N(E;\exp(2\pi r/N))}{N}
  =0.
\end{equation*}
\end{thm}
\begin{proof}
In this case
\begin{align*}
  g(j)
  &=4\sinh(\pi|r|j/N+\pi|r|)\sinh(\pi|r|j/N-\pi|r|)
  \\
  &=2\cosh(2\pi|r|)-2\cosh(2\pi|r|j/N)
  >0.
\end{align*}
Since $g(j)$ is monotonically decreasing and $g(0)=2\cosh(2\pi|r|)-2$,
$f(k)=\prod_{j=1}^{k}g(j)$ takes its maximum at
$g(\floor{\varphi(|r|)N/2\pi|r|})$
if $g(0)>1$.
Therefore
\begin{align*}
  &\lim_{N\to\infty}
  \frac{\log J_N{E;\exp(2\pi|r|/N)}}{N}
  \\
  &\quad=
  \lim_{N\to\infty}
  \frac{\log f(\floor{\varphi(|r|)N/2\pi|r|})}{N}
  \\
  &\quad=
  \lim_{N\to\infty}
  \sum_{j=1}^{\floor{\varphi(|r|)N/2\pi|r|}}
  \frac{\log 2\sinh(\pi|r|j/N+\pi|r|)}{N}
  \\
  &\qquad
  +
  \lim_{N\to\infty}
  \sum_{j=1}^{\floor{\varphi(|r|)N/2\pi|r|}}
  \frac{\log 2\sinh(\pi|r|j/N-\pi|r|)}{N}
  \\
  &\quad=
  \frac{1}{\pi|r|}\int_{\pi|r|}^{\varphi(|r|)/2+\pi|r|}\log 2\sinh x d\,x
  +
  \frac{1}{\pi|r|}\int_{-\pi|r|}^{\varphi(|r|)/2-\pi|r|}\log 2\sinh x d\,x
\end{align*}
and the result follows when $2\pi|r|>\arccosh(3/2)$.
\par
If $2\pi|r|\le\arccosh(3/2)$, then $g(j)<g(0)=2\cosh(2\pi|r|)-2\le1$ and so
\begin{align*}
  N>J_N(J;\exp(2\pi|r|/N))>g(0)-\varepsilon=2\cosh(2\pi r)-2-\varepsilon
\end{align*}
for any $\varepsilon>0$ if $N$ is sufficiently large.
Therefore
\begin{equation*}
  \lim_{N\to\infty}
  \frac{\log J_N(E;\exp(2\pi r/N))}{N}
  =0
\end{equation*}
as required.
\end{proof}
\begin{ack}
Part of this work was done when the author was visiting the International
Institute for Advanced Studies from 5th to 8th March, 2002 to attend the
workshop `Volume Conjecture and Its Related Topics' (partially supported by
the Research Institute for Mathematical Sciences, Kyoto University),
Warwick University to attend the workshop `Quantum Topology' from 18th to 22nd March, 2002,
and Universit{\'e} du Qu{\'e}bec {\`a} Montr{\'e}al to attend the workshop
`Knots in Montreal II' from 20th to 21st April, 2002.
Thanks are due to the institute/universities and to the organizers of the
workshops,
Stavros Garoufalidis, Colin Rourke, Steven Boyer, and Adam Sikora.
\par
The author thanks K.~Hikami for introducing the computer program PARI.
He is also grateful to S.~Gukov, A.~Mednykh, and A.~Vesnin for informing
their results.
\end{ack}
\appendix
\section{Some calculations}
\begin{lem}
For $|r-1|<1/3$, we have
\begin{equation*}
  \Lambda(\pi r+\theta(r)/2)-\Lambda(\pi r-\theta(r)/2)>0
\end{equation*}
\end{lem}
\begin{proof}
Put $V(r):=\Lambda(\pi r+\theta(r)/2)-\Lambda(\pi r-\theta(r)/2)$
Since
\begin{equation*}
  4\sin(\pi r+\theta(r)/2)\sin(\pi r-\theta(r)/2)
  =
  2\cos(\theta(r))-2\cos(2\pi r)
  =-1,
\end{equation*}
we have
\begin{align*}
  \frac{d\,V(r)}{d\,r}
  &=
  -\log\left|2\sin(\pi r+\theta(r)/2)\right|
  \left(\pi+\frac{1}{2}\frac{d\,\theta(r)}{d\,r}\right)
  \\
  &\quad
  +\log\left|2\sin(\pi r-\theta(r)/2)\right|
  \left(\pi-\frac{1}{2}\frac{d\,\theta(r)}{d\,r}\right)
  \\
  &=
  2\pi\log\left|2\sin(\pi r-\theta(r)/2)\right|.
\end{align*}
Since $\pi/6\le\pi r-\theta(r)/2<5\pi/6$ for $2/3<r<1$,
$\sin(\pi r-\theta(r)/2)>1/2$ and so $d\,V(r)/d\,r>0$.
Therefore $V(r)\ge0$ for $2/3<r<1$ since $V(2/3)=0$.
\par
Since $V(2-r)=V(r)$, the equality also holds for $1<r<4/3$.
\end{proof}
\begin{lem}
For $1/6<|r-1|<1/4$, we have
\begin{equation*}
  \Lambda(\pi r+\theta(r)/2)-\Lambda(\pi r-\theta(r)/2)
 +\Lambda(\pi r+\varphi(r)/2)-\Lambda(\pi r-\varphi(r)/2)>0
\end{equation*}
\end{lem}
\begin{proof}
Put $W(r):=\Lambda(\pi r+\theta(r)/2)-\Lambda(\pi r-\theta(r)/2)
 +\Lambda(\pi r+\varphi(r)/2)-\Lambda(\pi r-\varphi(r)/2)$.
As in the proof of the previous lemma, we have
\begin{align*}
  \frac{d\,W(r)}{d\,r}=
  2\log|2\sin(\pi r-\theta(r)/2)|
  +
  2\log|2\sin(\pi r-\varphi(r)/2)|.
\end{align*}
Since $7\pi/12<\pi r-\varphi(r)/2<5\pi/6$ for $3/4<r<5/6$,
$\sin(\pi r-\varphi(r)/2)>1/2$ and so $d\,W(r)/d\,r>0$.
Therefore $W(r)>0$ for $3/4<r<5/6$ since $V(3/4)=0$.
\par
Since $W(2-r)=W(r)$, the equality also holds for $7/6<r<5/4$.
\end{proof}
\bibliography{mrabbrev,hitoshi}

\providecommand{\bysame}{\leavevmode\hbox to3em{\hrulefill}\thinspace}
\providecommand{\MR}{\relax\ifhmode\unskip\space\fi MR }
\providecommand{\MRhref}[2]{%
  \href{http://www.ams.org/mathscinet-getitem?mr=#1}{#2}
}
\providecommand{\href}[2]{#2}
\begin{thebibliography}{10}

\bibitem{Gukov:03}
S.~Gukov, \emph{{Three-Dimensional Quantum Gravity, Chern-Simons Theory, and
  the A-Polynomial}}, HUTP-03/A003, ITEP-TH-50/02, arXiv:hep-th/0306165.

\bibitem{Kashaev:MODPLA95}
R.~M. Kashaev, \emph{A link invariant from quantum dilogarithm}, Modern Phys.
  Lett. A \textbf{10} (1995), no.~19, 1409--1418. \MR{96j:81060}

\bibitem{Kashaev:LETMP97}
\bysame, \emph{The hyperbolic volume of knots from the quantum dilogarithm},
  Lett. Math. Phys. \textbf{39} (1997), no.~3, 269--275. \MR{98b:57012}

\bibitem{Kashaev/Tirkkonen:ZAPNS200}
R.~M. Kashaev and O.~Tirkkonen, \emph{A proof of the volume conjecture on torus
  knots}, Zap. Nauchn. Sem. S.-Peterburg. Otdel. Mat. Inst. Steklov. (POMI)
  \textbf{269} (2000), no.~Vopr. Kvant. Teor. Polya i Stat. Fiz. 16, 262--268,
  370. \MR{1 805 865}

\bibitem{Murakami:4_1}
H.~Murakami, \emph{The asymptotic behavior of the colored {J}ones function of a
  knot and its volume}, Proceedings of `Art of Low Dimensional Topology VI'
  (T.~Kohno, ed.), January 2000, arXiv:math.GT/0004036.

\bibitem{Murakami:SURIK00}
\bysame, \emph{Optimistic calculations about the
  {W}itten-{R}eshetikhin-{T}uraev invariants of closed three-manifolds obtained
  from the figure-eight knot by integral {D}ehn surgeries},
  S\=urikaisekikenky\=usho K\=oky\=uroku (2000), no.~1172, 70--79. \MR{1 805
  729}

\bibitem{Murakami:Nagoya99}
\bysame, \emph{Kashaev's invariant and the volume of a hyperbolic knot after
  {Y}. {Y}okota}, Physics and combinatorics 1999 (Nagoya), World Sci.
  Publishing, River Edge, NJ, 2001, pp.~244--272. \MR{2002k:57031}

\bibitem{Murakami:SURIK02}
\bysame, \emph{Mahler measure of the colored {J}ones polynomial and the volume
  conjecture}, S\=urikaisekikenky\=usho K\=oky\=uroku (2002), no.~1279, 86--99.
  \MR{1 953 834}

\bibitem{Murakami/Murakami:ACTAM101}
H.~Murakami and J.~Murakami, \emph{The colored {J}ones polynomials and the
  simplicial volume of a knot}, Acta Math. \textbf{186} (2001), no.~1, 85--104.
  \MR{1 828 373}

\bibitem{Murakami/Murakami/Okamoto/Takata/Yokota:EXPMA02}
H.~Murakami, J.~Murakami, M.~Okamoto, T.~Takata, and Y.~Yokota, \emph{Kashaev's
  conjecture and the {C}hern-{S}imons invariants of knots and links},
  Experiment. Math. \textbf{11} (2002), no.~3, 427--435. \MR{1 959 752}

\bibitem{Vesnin/Mednykh:1995}
A.~Yu. Vesnin and A.~D. Mednykh, \emph{Hyperbolic volumes of {F}ibonacci
  manifolds}, Sibirsk. Mat. Zh. \textbf{36} (1995), no.~2, 266--277, i.
  \MR{97c:57019}

\bibitem{Yokota:volume00}
Y.~Yokota, \emph{{On the volume conjecture for hyperbolic knots}},
  arXiv:math.QA/0009165.

\bibitem{Yokota:Topology_Symposium2000}
\bysame, \emph{On the volume conjecture for hyperbolic knots}, Proceedings of
  the 47th Topology Symposium (Inamori Hall, Kagoshima University), July 2000,
  pp.~38--44.

\bibitem{Yokota:Murasugi70}
\bysame, \emph{On the volume conjecture of hyperbolic knots}, Knot Theory --
  dedicated to Professor Kunio Murasugi for his 70th birthday (M.~Sakuma, ed.),
  March 2000, pp.~362--367.

\bibitem{Yokota:GTM02}
\bysame, \emph{{On the potential functions for the hyperbolic structures of a
  knot complement}}, Geom. Topol. Monogr. \textbf{4} (2002), 303--311.

\bibitem{Yokota:INTIS03}
\bysame, \emph{From the {J}ones polynomial to the {$A$}-polynomial of
  hyperbolic knots}, Interdiscip. Inform. Sci. \textbf{9} (2003), 11--21.

\end{thebibliography}
\bibliographystyle{amsplain}
\end{document}